\documentclass[12pt,leqno]{article}
\newcommand{\EE}{{\rm I\kern-2pt E}}
\newcommand{\RR}{{\rm I\kern-2pt R}}
\newcommand{\DD}{{\rm I\kern-2pt D}}
\newcommand{\PP}{{\rm I\kern-2pt P}}
\newcommand{\NN}{{\rm I\kern-2pt N}}
\newcommand{\dd}{{\rm \kern 3pt I\kern-9pt d}}

\title{\large THEOREME DE DONSKER ET FORMES DE DIRICHLET}
\author{\sc Nicolas Bouleau}
\date{\it Ecole des Ponts, Paris}

\begin{document}
\maketitle

\noindent{\bf Abstract.} We use the language of errors to handle local Dirichlet forms with squared field operator (cf [2]). Let us consider, under the hypotheses of Donsker theorem, 
a random walk converging weakly to a Brownian motion. If, in addition, the random walk is  supposed to be erroneous, the convergence occurs in the sense
of Dirichlet forms and induces the Ornstein-Uhlenbeck structure on the Wiener space. This quite natural result uses an extension of Donsker theorem
to functions with quadratic growth. As an application we prove an invariance principle for the gradient of the maximum of the Brownian path computed
by Nualart and Vives.\\

\noindent{\bf R\'esum\'e.} Nous employons le langage des erreurs pour manier les formes de Dirichlet locales avec carr\'e du champ (cf [2]). Consid\'erant
une promenade al\'eatoire convergeant en loi vers un mouvement brownien sous les hypoth\`eses du th\'eor\`eme de Donsker, nous montrons que si la promenade
est suppos\'ee de plus erron\'ee, la convergence a lieu au sens des formes de Dirichlet et induit la structure d'Onstein-Uhlenbeck sur l'espace de Wiener. Ce 
r\'esultat bien naturel n\'ecessite l'extension du th\'eor\`eme de Donsker aux fonctions \`a croissance quadratique. A titre d'application nous en d\'eduisons 
un principe d'invariance pour le gradient du maximum de la courbe brownienne calcul\'e par Nualart et Vives.\\

Mots cl\'es : promenade al\'eatoire, mouvement brownien, gradient, forme de Dirichlet, erreur.

Keywords : random walk, Brownian motion, Dirichlet form, error\\

\noindent{1. INTRODUCTION}\\

Le calcul d'erreur fond\'e sur les formes de Dirichlet est  inspir\'e des id\'ees de Gauss sur la propagation des erreurs en les formulant
avec des formes de Dirichlet ce qui leur donne la puissance de s'appliquer aux espaces fonctionnels rencontr\'es en mod\'elisation stochastique :
 espace de Wiener, 
 de Poisson et de Monte Carlo. L'approche intuitive et les fondements math\'ematiques sont expos\'es en [2]. Ce formalisme rend compte de la 
propagation \`a travers les calculs des variances, des co-variances et des biais des erreurs suppos\'ees infiniment petites.

De la m\^eme fa\c{c}on que le calcul des probabilit\'es a pu se d\'evelopper sans que la notion de hasard f\^ut compl\^etement \'elucid\'ee, 
le calcul d'erreur fond\'e sur les formes de Dirichlet n'explicite pas la notion d'erreur elle-m\^eme et ne prend en compte que les notions d\'eriv\'ees de variance et de biais
qui sont axiomatis\'ees. L'explicitation de cette notion d'erreur est un programme th\'eorique int\'eressant mais n'est pas un pr\'ealable aux nombreuses
 applications des calculs d'erreurs et de sensibilit\'e (cf. [4], [5] et [2]).

En revanche la question du choix des hypoth\`eses \`a prendre en compte sur les erreurs lorsqu'on proc\`ede \`a une \'etude de sensibilit\'e d'un mod\`ele 
est une question concr\`ete importante dans la mesure o\`u ce langage plus fin que d'ordinaire sur les erreurs fait appara\^{\i}tre la n\'ecessit\'e
d'hypoth\`eses a priori sur les corr\'elations ou non-corr\'elations des erreurs sur les param\`etres scalaires ou fonctionnels du mod\`ele \'etudi\'e.
La connexion avec les statistiques par l'information de Fisher et sa robustesse par changement de variables est la r\'eponse g\'en\'erale \`a cette question
(cf. [6]). Elle peut \^etre compl\`et\'ee par l'\'etude des extensions en termes de calcul d'erreur des grands th\'eor\`emes limites de la th\'eorie des 
probabilit\'es tels que le th\'eor\`eme de limite centrale ou le th\'eor\`eme du logarithme it\'er\'e (cf. [7]).

Nous \'etudions ici la question bien naturelle de l'extension du th\'eor\`eme de Donsker concernant la limite faible d'une promenade al\'eatoire vers un
mouvement brownien. Elle peut se formuler ainsi : \'etant donn\'e une suite de variables al\'eatoires ind\'ependantes
 \'equidistribu\'ees centr\'ees, suppos\'ees en outre erron\'ees, les erreurs \'etant stationnaires et non corr\'el\'ees, est-ce que l'approximation usuelle
affine par morceaux converge vers le mouvement brownien au sens de la forme de Dirichlet qui d\'ecrit les erreurs et si oui, quelle structure d'erreur cela
induit-il sur l'espace de Wiener ?

La r\'eponse est positive et la structure d'erreur obtenue est la structure d'Orn\-stein-Uhlenbeck. Ce r\'esultat tr\`es naturel n'avait pas 
\'et\'e publi\'e jusqu'ici, sa d\'emons\-tration n\'ecessite une am\'elioration strictement proba\-biliste du th\'eor\`eme de Donsker  aux fonctions
\`a croissance quadratique. Cette extension, plus d\'elicate que dans le cas du th\'eor\`eme de limite centrale est la principale difficult\'e du pr\'esent travail.

Nous en tirons comme cons\'equence,  une formule explicite concernant la limite de la forme de Dirichlet sur la norme uniforme
qui utilise le beau r\'esultat de Nualart et Vives [11] sur le gradient du maximum de la trajectoire brownienne sur [0,1].\\

\noindent{2. DEFINITIONS ET NOTATIONS}\\

Une structure d'erreur est un terme $S=(\Omega,\mathcal{A},\PP,\DD,\Gamma)$  o\`u $(\Omega,\mathcal{A},\PP)$ est
 un espace de probabilit\'e, $\DD$ un sous-espace dense de   $L^2(\PP)$ et  $\Gamma$ un op\'erateur bilin\'eaire
 sym\'etrique positif de      $\DD\times\DD$ dans $L^1(\PP)$ v\'erifiant 

1) le calcul fonctionnel de classe  $C^1 \cap Lip$

 i.e. si  $U=(U_1,\ldots,U_m)\in \DD^ m $, $V=(V_1,\ldots,V_n)\in \DD^n $,  

et si $F$ et $G$ sont de classe $ C^1$ et lipschitziennes
 de $\RR^m$ [resp.$\RR^n$] dans $\RR$, 

alors $F(U_1,\ldots,U_m)\in\DD$ et $G(V_1,\ldots,V_n)\in \DD$ et 

$\Gamma[F(U_1,\ldots,U_m),G(V_1,\ldots,V_n)]=\sum\limits_{i,j} F_i'(U) G_j'(V)\Gamma[U_i,V_j]\quad \PP{\mbox{-p.s.}}$ 

2) $1\in\DD$,  la forme    $\mathcal{E}[F,G]=\frac{1}{2}\int \Gamma[F,G] d\PP $ est ferm\'ee

i.e.  $\DD$  est complet pour la norme $\parallel. \parallel_{\mathcal{E}}=(\parallel.\parallel^2_{L^2(\PP)}+\mathcal{E}[.])^{\frac{1}{2}}$.\\

Avec les notations ci-dessus $\mathcal{E}$ est une forme de Dirichlet locale ayant pour op\'erateur carr\'e du champ $\Gamma$. 
On notera $\Gamma[F]$ pour $\Gamma[F,F]$ et $\mathcal{E}[F]$ pour $\mathcal{E}[F,F]$.\\

\noindent{\it Exemple.} Si l'on note $\mathcal{B}(\RR)$ la tribu bor\'elienne sur $\RR$, $\nu=\mathcal{N}(0,1)$ la 
loi normale centr\'ee r\'eduite et $H^1(\nu)$ l'espace de Sobolev associ\'e, alors 
$$(\RR,\mathcal{B}(\RR),\nu,H^1(\nu), \Gamma[u]=u^{'2}) $$ est une structure d'erreur appel\'ee 
structure d'erreur d'Ornstein-Uhlenbeck sur $\RR$.\\

Les structures d'erreur ont la propri\'et\'e de se transporter simplement par image et de permettre les op\'erations 
de produits y compris infini-d\'enombrables, cf. [2], [7].\\

\noindent{\it Convergence en loi de Dirichlet.}

Soit une structure d'erreur $S=(\Omega,\mathcal{A},\PP,\DD,\Gamma)$, soit ${\cal E}$ la forme de Dirichlet associ\'ee.

Soit $W$ un espace vectoriel norm\'e muni de sa tribu bor\'elienne ${\cal B}(W)={\cal W}$.

On se donne une famille de variables al\'eatoires  $(U_n)_{n\in\NN}$ d\'efinies sur $(\Omega, {\cal A})$  \`a valeurs $(W,\mathcal{W})$.

 On introduit une notion de convergence adapt\'ee pour les structures d'erreur de la convergence en loi des variables al\'eatoires.\\ 

\noindent{\bf D\'efinition 1}. {\it On dit que  $(U_n)_{n\in\NN}$ converge en loi de Dirichlet 
 s'il existe une structure d'erreur sur $(W, {\cal W})$  soit  $\Sigma=(W, {\cal W},m,\DD_0,\Gamma_0)$ telle que  :\\

i) $(U_n)_*\PP\rightarrow m$  \'etroitement  

i.e.  $\forall f:\Omega\mapsto \RR$ 
continue born\'ee $\EE[f(U_n)]\longrightarrow \int_{\cal W} f(w) dm(w)$,\\

ii) si $ F\in {\cal C}^1\cap Lip(W,\RR)$ alors $F\in\DD_0$ et $F(U_n)\in\DD\;\forall n$ et

$\mathcal{E}[F(U_n)]\longrightarrow\mathcal{E}_0[F]\quad{\mbox{quand}}\quad n\rightarrow\infty.$
 o\`u ${\cal E}_0$ est la forme associ\'ee \`a $\Sigma$.}\\

\noindent{\it Remarque 1}. Sous les hypoth\`eses de la d\'efinition 1, les $U_n$ transportent la structure $S$ sur $(W,{\cal W})$:

Si on d\'efinit
$$\begin{array}{rcl}
\PP_{U_n}&=&(U_n)_*\PP\quad {\mbox{(loi de }} U_n{\mbox{)}}\\
\DD_{U_n}&=&\{\varphi\in L^2(\PP_{U_n})\,:\,\varphi(U_n)\in\DD\}\\
\Gamma_{U_n}[\varphi](w)&=&\EE[\Gamma[\varphi(U_n)]|U_n=w]
\end{array}
$$
\noindent le terme
$$S_{U_n}= (W,{\cal W},\PP_{U_n},\DD_{U_n}, \Gamma_{U_n})$$ est une structure d'erreur, $\DD_{U_n}$ contient les 
fonctions $C^1\cap Lip(W,\RR)$, $\PP_{U_n}$ converge \'etroitement vers $m$ sur $(W,{\cal W})$ et
${\cal E}_{U_n}[F]=\frac{1}{2}\int \Gamma_{U_n}[F]\,d\PP_{U_n}\rightarrow\frac{1}{2} \int \Gamma_0[F]\,dm$
pour toute $F\in C^1\cap Lip(W,\RR)$. Il est naturel d'appeler la structure $S_{U_n}$ la loi de Dirichlet de $U_n$.\\

\noindent{\it Remarque 2}. Si de plus il existe une variable al\'eatoire $V\in\DD_0$ telle que

i) $(U_n)_*\PP\rightarrow V_*m$ \'etroitement

ii) $\forall F\in C^1\cap Lip\quad {\cal E}[F(U_n)]\rightarrow {\cal E}_0[F(V)]$

\noindent nous dirons que les $U_n$ convergent en loi de Dirichlet vers $V$.\\

\noindent 3. CONVERGENCE D'UNE PROMENADE ALEATOIRE ERRONEE.	\\

Rappelons le r\'esultat classique de Donsker [8] concernant la convergence d'une promenade al\'eatoire. Soient $U_n,\,n\geq 1$, une suite de variables
al\'eatoires i.i.d. de carr\'e int\'egrable, de variance $\sigma^2$, centr\'ees. On interpole la promenade al\'eatoire $\sum_{k=1}^nU_k$ de fa\c{c}on affine par morceaux
en consid\'erant le processus
$$X_n(t)=\frac{1}{\sqrt{n}}\left(\sum_{k=1}^{[nt]}U_k+(nt-[nt])U_{[nt]+1}\right)$$
pour $t\in[0,1]$, o\`u $[x]$ d\'esigne la partie enti\`ere de $x$.

L'espace $W=C([0,1])$ \'etant muni de la norme uniforme, les variables $X_n$ \`a valeur $W$ convergent en loi vers un mouvement brownien centr\'e de 
variance $\sigma^2 t$.

Il en r\'esulte que si $\Phi$ est une fonctionnelle Riemann-int\'egrable pour la mesure de Wiener et born\'ee
$$\EE[\Phi(X_n)]\rightarrow\EE[\Phi(B)],$$ o\`u $B$ est un mouvement brownien centr\'e de variance $\sigma^2 t$.\\

Supposons maintenant que les $U_n$ soient erron\'ees, en conservant les hypoth\`eses d'ind\'ependance et d'\'equi-distribution pour les $U_n$ et leurs erreurs.
Autrement dit, consid\'erons que les $U_n$ sont les applications coordonn\'ees d'une structure d'erreur produit
$$S=(\Omega, {\cal A}, \PP, \DD,\Gamma)=(\RR, {\cal B}(\RR), \mu, {\bf d},\gamma)^{\NN^\ast}$$
la structure $(\RR, {\cal B}(\RR), \mu, {\bf d},\gamma)$ \'etant telle que l'identit\'e $j$ soit dans $L^2(\mu)$ centr\'ee et dans {\bf d}. Ainsi les $U_n$
sont i.i.d., de loi $\mu$, de variance $\sigma^2=\mu(j^2)$, v\'erifient $U_n\in\DD$ et
\begin{eqnarray}
\left\{
\begin{array}{rcl}
\Gamma[U_n]&=&(\gamma[j])(U_n)\\
\Gamma[U_m,U_n]&=&0\;{\mbox{ si }} m =\!\!\!\!\!/\;\, n
\end{array}
\right.
\end{eqnarray}
Les v.a. $\Gamma[U_n]$ sont dans $L^1(\PP)$, ind\'ependantes et de m\^eme loi.

Pour $t$ fix\'e la v. a. 
$$X_n(t)=\frac{1}{\sqrt{n}}\left(\sum_{k=1}^{[nt]}U_k+(nt-[nt])U_{[nt]+1}\right)$$
est dans $\DD$, et par (1) 
\begin{eqnarray}
\Gamma[X_n(s),X_n(t)]=\frac{1}{n}\left[\sum_{k=1}^{[nt]\wedge[ns]}\Gamma[U_k]+\alpha(n,s,t)\right] 
\end{eqnarray}
avec
$$
\begin{array}{rl}
\alpha(n,s,t)=&\left((ns-[ns])1_{\{[ns]<[nt]\}}+(nt-[nt])1_{\{[ns]>[nt]\}}\right.\\
&+\left.(ns-[ns])(nt-[nt])1_{\{[ns]=[nt]\}}\right)\Gamma[U_{[ns]\wedge[nt]+1}]
\end{array}
$$
Il d\'ecoule de la loi forte des grands nombres que
$$
\frac{1}{n}\sum_{k=1}^{[nt]\wedge[ns]}\Gamma[U_k]\rightarrow (s\wedge t)\EE[\Gamma[U_1]]\quad \PP{\mbox{-p.s. et dans }} L^1(\PP).$$
Par ailleurs 
$$\frac{|\alpha(n,s,t)|}{n}\rightarrow 0 \quad \PP{\mbox{-p.s. et dans }} L^1(\PP).$$
Ainsi $\Gamma[X_n(s),X_n(t)]\rightarrow(s\wedge t)c \;\; \PP{\mbox{-p.s. et dans }} L^1(\PP)$ o\`u $c$ est
 la constante $\EE \Gamma[U_1]=\int \gamma[j](x)d\mu(x)$.

Nous d\'eduisons de ce calcul la convergence des lois de Dirichlet marginales d'ordre fini vers les marginales correspondantes du mouvement brownien
muni de la structure d'Ornstein-Uhlenbeck: soit $W=C([0,1])$ muni de sa tribu bor\'elienne ${\cal W}$ et $m$ la mesure de Wiener telle que la coordonn\'ees d'indice $t$ soit centr\'ee 
de variance $\sigma^2 t$, soit $\DD_0$ le domaine de la forme d'Ornstein-Uhlenbeck et $\Gamma_0$ l'op\'erateur quadratique associ\'e caract\'eris\'e par son
action sur le premier chaos (cf. [2] chapitre VI \S 2 et [7])
$$\forall h\in L^2([0,1])\quad\int_0^1hdB\in\DD_0 \;{\mbox{ et }}\;\Gamma_0[\int_0^1hdB]=c\int h^2dt.$$

\noindent{\bf Proposition 1.} {\it Soient $t_1,\ldots,t_p\in[0,1]$,  les variables al\'eatoires 
$(X_n(t_1),\ldots, X_n(t_p))$ convergent en loi de Dirichlet vers $(B(t_1),\ldots, B(t_p))$ o\`u $B$ est un mouvement brownien centr\'e de variance $\sigma^2 t$
muni de la structure d'Ornstein-Uhlenbeck

\noindent$(W,{\cal W},m,\DD_0,\Gamma_0)$.}\\

\noindent{\it D\'emonstration}. Il faut montrer que si $f\in C^1\cap Lip$ 
$$\int\Gamma[f(X_n(t_1),\ldots, X_n(t_p))]\,d\PP\rightarrow\int \Gamma_0[f(B(t_1),\ldots, B(t_p))]\,dm.$$
Par majoration de la fonction $\alpha(n,t_i,t_j)$ et par le calcul fonctionnel on est ramen\'e \`a \'etudier la convergence de l'expression
\begin{eqnarray}
\EE[f_i^\prime((X_n(t_1),\ldots, X_n(t_p))f_j^\prime((X_n(t_1),\ldots, X_n(t_p))\frac{1}{n}\sum_{k=1}^{[nt_i]\wedge[nt_j]}\Gamma[U_k]]
\end{eqnarray}
et pour cela d'\'etudier la convergence de 
$$\EE[e^{i(u_1X_n(t_1)+\cdots+u_pX_n(t_p))}\Gamma[U_k]]$$ pour $k$ fix\'e. Or, compte tenu de ce que $\Gamma[U_k]=(\gamma[j])(U_k)$,
par un argument classique (similaire ˆ celui de la proposition 3 ci-dessous), cette expression converge
vers
$$\EE[e^{i(u_1B(t_1)+\cdots+u_pB(t_p))}]c.$$ D'o\`u finalement (3) converge vers
$$\EE[f_i^\prime((B(t_1),\ldots, B(t_p))f_j^\prime((B(t_1),\ldots, B(t_p))]c(t_i\wedge t_j)$$ ce qui d\'emontre la proposition.\\

Ces r\'esultats sur les marginales finies posent naturellement la question de l'extension suivante du th\`eor\`eme de Donsker :\\

\noindent{\bf Th\'eor\`eme 1.} {\it Les variables $X_n$ convergent en loi de Dirichlet vers la structure d'ornstein-Uhlenbeck sur l'espace de Wiener 
$(W,{\cal W},m,\DD_0,\Gamma_0)$.}\\

Nous donnerons deux d\'emonstrations de ce th\'eor\`eme. La premi\`ere plus \'el\'emen\-taire n\'ecessite l'hypoth\`ese suppl\'ementaire que la fonction $\gamma[j]$
est dans  $L^p(\mu)$ pour un $p>1$. Elle fait comprendre la difficult\'e surmont\'ee par la
 seconde d\'emonstration qui utilise un renforcement du th\'eor\`eme de Donsker probabiliste. Nous aurons besoin de quelques lemmes et notations.

\noindent{\bf Lemme 1.} {\it Si $F\in C^1\cap Lip(W,\RR)$ o\`u $W$ est muni de la norme uniforme,
$$F(x+h)=F(x)+<F'(x),h>+\| h\| \varepsilon_x(h)$$ o\`u $\varepsilon_x(h)$ est born\'ee (en $x$ et $h$) et $\varepsilon_x(h)\rightarrow 0$ quand $h\rightarrow 0$
dans $W$, et o\`u $x\mapsto F'(x)$ est continue born\'ee de $W$ dans l'espace de Banach des mesures de Radon sur $[0,1]$.}

\noindent{\it D\'emonstration}. Il r\'esulte en effet des hypoth\`eses que $|<F'(x),u>|\leq K$ pour tout $u$ unitaire dans $W$ o\`u $K$ est la constante de Lipschitz de $F$ 
ce qui donne le 
r\'esultat.\\

Il sera commode d'utiliser l'op\'erateur $(.)^\#$ qui est un gradient particulier  cons\-truit avec une copie de l'espace initial (cf. [2] p80).

Soit $(\hat{\Omega}, \hat{\cal A},\hat{\PP})$ une copie de $(\Omega, {\cal A}, \PP)$  et $\widehat{U_n}$ les coordonn\'ees de $\hat{\Omega}$. 
Choisissant un op\'erateur di\`ese pour la structure $(\RR, {\cal B}(\RR), \mu, {\bf d},\gamma)$, nous en d\'eduisons un op\'erateur di\`ese pour la structure 
produit (cf. [2] p80 remarque) en posant $U_n^{\#}=j^{\#}(U_n,\widehat{U_n})$.
 Maintenant  pour d\'efinir l'op\'erateur $(.)^{\#}$ de
$\DD$ dans $L^2 (\Omega\times\hat{\Omega}, \PP\times\hat{\PP})$ il n'est que de poser si $H=h(U_1,\ldots,U_k,\ldots)\in \DD$ 
$$H^{\#}=\sum_i h'_i(U_1,\ldots,U_n,\ldots)U_i^{\#}.$$
On a alors $$\hat{\EE}[(H^\#)^2]=\Gamma[H]\qquad\forall H\in\DD$$ d'o\`u il r\'esulte (cf. [2]) que $\forall \varphi\in C^1\cap Lip, \;\forall H_1,\ldots,H_p\in \DD$
$$(\varphi(H_1,\ldots,H_p))^\#=\sum_i\varphi'_i(H_1,\ldots, H_q)H_i^\#.$$
De m\^eme sur l'espace de Wiener, nous consid\'erons une copie $(\hat{W},\hat{\cal W},\hat{m})$ et l'op\'erateur $(.)^\#$ de $\DD_0$ dans $L^2(W\times{\cal W},m\times\hat{m})$
qui v\'erifie (cf. [2] chap VI \S 2) $(B_t)^\#=\frac{\sqrt{c}}{\sigma}\hat{B}_t$ et $\forall H\in\DD_0\quad\hat{\EE}[(H^\#)^2]=\Gamma_0[H].$

Nous avons alors

\noindent{\bf Lemme 2.} {\it Soit $F\in C^1\cap Lip(W)$, on a $F(X_n)\in\DD$ et
$$(F(X_n))^\#=\int_{[0,1]}(X_n(s))^\#\,F'(X_n)(ds)$$
et 
$$\Gamma[F(X_n)]=\int_{[0,1]}\int_{[0,1]}\Gamma[X_n(s),X_n(t)]F'(X_n)(ds)F'(X_n)(dt).$$
De m\^eme $F(B)\in\DD_0$ et 
$(F(B))^\#=\int_{[0,1]}B^\#(s)\,F'(B)(ds)$ et
$$
\Gamma[F(B)]=\int_{[0,1]}\int_{[0,1]}s\wedge t \,F'(B)(ds)F'(B)(dt)
=\int_0^1<F'(B),1_{[u,1]}>^2\,du.
$$}
{\it D\'emonstration.}  Les formules sont ais\'ees \`a \'etablir lorsque $F(x)$ ne d\'epend que d'un nombre fini des valeurs prises par $x$. Puis pour $F\in
C^1\cap Lip$ quelconque, soit $P_k(x)$ l'approximation de $x$ affine par morceaux de pas $\frac{1}{k}$ et posons $F_k(x) = F\circ P_k(x)$. $F_k$ est de 
classe $C^1\cap Lip$ et on a $(F_k)'=P_k^\ast F'P_k$ o\`u $P_k^\ast$ est l'op\'erateur adjoint de $P_k$. Nous avons
$$(F_k(X_n))^\#=\int_{[0,1]}(X_n(s))^\#\,F'_k(X_n)(ds).$$
Compte tenu de ce que $s\mapsto (X_n(s))^\#$ est continue comme \'etant l'approximation affine par morceau de la promenade des $U^\#_n$, nous avons
$$(F_k(X_n))^\#=<P_k(X_n(.))^\#,F'(P_k(X_n))>.$$
D'o\`u il r\'esulte que $F_k(X_n))^\#$ converge vers $<(X_n(.))^\#,F'(X_n)>$ en restant domin\'e en module 
par $\max_s |(X_n(s))^\#|\max_x\|F'(x)\|$ donc dans $L^2(\PP\times\hat{\PP})$.
La premi\`ere formule r\'esulte alors du fait que $(.)^\#$ est un op\'erateur ferm\'e. Les autres formules en d\'ecoulent ou 
se d\'emontrent de fa\c{c}on analogue.\\

La premi\`ere d\'emonstration du th\'eor\`eme s'engage maintenant naturellement:

\noindent{\it Premi\`ere d\'emonstration}

Du lemme pr\'ec\'edent et de la formule (2) nous tirons 
$$
\begin{array}{rl}
\Gamma[F(X_n)]=&\int\int\left(\frac{1}{n}\sum_{k=1}^{[ns]\wedge[nt]}\Gamma[U_k]\right)F'(X_n)(ds)F'(X_n)(dt)\\
&+ \int\int\alpha(n,s,t)F'(X_n)(ds)F'(X_n)(dt)\\
=&(A)+(B).
\end{array}
$$
Le second terme peut \^etre major\'e ainsi
$$|(B)|\leq\frac{1}{n}\sup_{k\leq n}\Gamma[U_k]\|F'(X_n)\|^2$$ o\`u $\|F'(X_n)\|$ est la masse totale de la mesure $F'(X_n)$. Il r\'esulte alors du lemme suivant 
que 
$\EE[|(B)|]\rightarrow 0$ quand $n$ tend vers l'infini.

\noindent{\bf Lemme 3.} {\it Si les $Y_k$ sont i.i.d. dans $L^1$ et positives, 
$\lim_n\EE[\frac{1}{n}\sup_{k\leq n} Y_k]= 0. $ }

\noindent{\it Preuve du lemme.} On a 
$$\EE[\frac{1}{n}\sup_{k\leq n} Y_k]=\int_0^\infty \frac{1-((1-\PP(Y_1>a))^n}{n}da$$ et 
$$\frac{1-((1-\PP(Y_1>a))^n}{n}\leq \PP(Y_1>a)$$ qui est int\'egrable puisque $Y_1\in L^1$ donc le th\'eor\`eme de convergence domin\'ee s'applique et 
donne le r\'esultat.

En ce qui concerne le premier terme $(A)$ posons $ V_k=\sum_{i=1}^k(\Gamma[U_i]-\EE\Gamma[U_i])$. 
Supposant $\Gamma[U_1]\in L^p$ pour un $p>1$, nous avons
par l'in\'egalit\'e de Doob ([10] p. 68) appliqu\'ee \`a la martingale $V_k$
$$\EE[\frac{1}{n}\max_{1\leq k\leq n}|V_k|]\leq \frac{p}{p-1}\frac{1}{n}\max_{1\leq k\leq n} k\|\frac{1}{k}\sum_{i=1}^k(\Gamma[U_i]-\EE\Gamma[U_i])\|_p$$
Le second membre est de la forme $\frac{1}{n}\max_{k\leq n} k\varepsilon(k)$ avec $\varepsilon(k)\rightarrow 0$ donc tend vers z\'ero quand $n\rightarrow\infty$.

Il r\'esulte de cette majoration que $\EE[(A)]$ a m\^eme limite que
$$\EE\int\int\left(\frac{1}{n}\sum_{k=1}^{[ns]\wedge[nt]}\EE[\Gamma[U_k]]\right)F'(X_n)(ds)F'(X_n)(dt)$$ ce qui vaut puisque les $\Gamma[U_k]$ sont
i.i.d.
$$c\EE\int_0^1<F'(X_n), 1_{[\frac{[nu]}{n},1]}>^2\,du$$
D'o\`u par le th\'eor\`eme de Donsker, l'application $x\mapsto \int<F'(x),1_{[u,1]}>^2du$ \'etant continue born\'ee, on a finalement 
$$\EE\Gamma[F(X_n)]\rightarrow c\EE\Gamma_0[F(B)]\hspace{3cm}{\mbox{C.Q.F.D.}}$$

Pour lever l'hypoth\`ese $\Gamma[U_1]\in L^p$ pour un $p>1$, nous abordons la question diff\'eremment.
De 
$$X_n(t)=\frac{1}{\sqrt{n}}\left(\sum_{k=1}^{[nt]}U_k+(nt-[nt])U_{[nt]+1}\right)$$
nous tirons
$$X_n^\#(t)=\frac{1}{\sqrt{n}}\left(\sum_{k=1}^{[nt]}U_k^\#+(nt-[nt])U_{[nt]+1}^\#\right)$$
de telle sorte que par le th\'eor\`eme de Donsker appliqu\'e aux couples $(U_k,U_k^\#)$ qui sont i.i.d. on a si $G$ est continue
 born\'ee de $W\times W$ dans $\RR$
$$\EE\hat{\EE}[G(X_n,X_n^\#)]\rightarrow\EE\hat{\EE}[G(B,B^\#)].$$
Pour  d\'emontrer le th\'eor\`eme en appliquant cette id\'ee \`a $F(X_n)^\#=\int X_n^\#(s)F'(X_n)(ds)$ il faudrait disposer 
du th\'eor\`eme de Donsker
non seulement pour les fonctions born\'ees mais pour les $G$ telles que
$|G(x)|\leq K_1\|x\|^2+K_2$. C'est ce que nous \'etablissons ici:\\

\noindent{\bf Th\'eor\`eme 2.} {\it Soient $X_n(t)$ comme dans le th\'eor\`eme de Donsker, et $B(t)$ un mouvement brownien de variance $\sigma^2 t$, 
alors $$\EE[\Phi(X_n)]\rightarrow\EE[\Phi(B)]$$ pour toute $\Phi$ continue de $W$ dans $\RR$ telle que $|\Phi(x)|\leq K_1\|x\|^2+K_2$.}

\noindent{\it D\'emonstration. } Posons
$$Z_n=\max_{t}|X_n(t)]=\frac{1}{\sqrt{n}}\max_{1\leq k\leq n}|\sum_{j=1}^k U_j|.$$
a) Il suffit de montrer que les v. a. $Z_n^2$ sont uniform\'ement int\'egrables.

En effet, $\varepsilon>0$ \'etant donn\'e, cette uniforme int\'egrabilit\'e entra\^{\i}ne que l'on peut trouver $a>0$ tel que 
$$|\EE[(\Phi(B)\wedge a)\vee(-a)]-\EE[\Phi(B)]|\leq \varepsilon/3$$ et que $\forall n$
$$|\EE[(\Phi(X_n)\wedge a)\vee(-a)]-\EE\Phi(X_n)|\leq\EE[|\Phi(X_n)|1_{|\Phi(X_n)|>a}]\leq \varepsilon/3.$$
Choisissant alors, d'apr\`es le th\'eor\`eme de Donsker, $n$ assez grand pour que 
$$|\EE[(\Phi(X_n)\wedge a)\vee(-a)]-\EE[(\Phi(B)\wedge a)\vee(-a)]|\leq \varepsilon/3$$ on a $|\EE\Phi(X_n)-\EE\Phi(B)|\leq \varepsilon.$

\noindent b) Pour montrer que les v. a. $Z_n^2$ sont uniform\'ement int\'egrables, nous posons $S_n=\sum_{i=1}^nU_i$ et nous utilisons la majoration suivante
([1]  p.69)
$$\PP\{\max_{i\leq n}|S_i|\geq \lambda\sigma\sqrt{n}\}\,\leq\,2\PP\{\frac{|S_n|}{\sigma\sqrt{n}}\geq \frac{\lambda}{2}\}\quad{\mbox{ si }}\lambda\geq 2\sqrt{2}$$
d'o\`u
$$\PP\{Z_n^2\geq \alpha\}\,\leq \,2\PP\{\frac{|S_n|}{\sigma\sqrt{n}}\geq \frac{\sqrt{\alpha}}{2\sigma}\}\quad{\mbox{ si }}\alpha\geq 8\sigma^2.$$
De 
$$\EE[Z_n^21_{Z_n^2\geq \alpha}]=\alpha\PP\{Z_n^2\geq\alpha\}+\int_\alpha^\infty\PP\{Z_n^2\geq t\}dt$$ nous d\'eduisons
$$\EE[Z_n^21_{Z_n^2\geq \alpha}]\leq 2\alpha\PP\{\frac{|S_n|}{\sigma\sqrt{n}}\geq \frac{\sqrt{\alpha}}{2\sigma}\}
+2\EE[(4\frac{S_n^2}{n}-\alpha)^+].$$ Il r\'esulte alors du th\'eor\`eme de limite centrale, et du fait que les $\frac{S_n^2}{n}$ sont uniform\'ement int\'egrables, 
que si $\alpha\geq 8\sigma^2$, 
\begin{eqnarray}
\limsup_{n}\EE[Z_n^21_{Z_n^2\geq\alpha}]\leq2\alpha\PP\{|N|\geq \frac{\sqrt{\alpha}}{2\sigma}\}+2\EE(4N^2-\alpha)^+
\end{eqnarray}
o\`u $N$ est une variable normale r\'eduite. Donc 
$$\lim_{\alpha\uparrow\infty}\limsup_{n}\EE[Z_n^21_{Z_n^2\geq\alpha}]=0$$ et ceci entra\^{\i}ne l'uniforme int\'egrabilit\'e des $Z_n^2$.\hspace{3cm}C.Q.F.D.\\

\noindent{\it Revenons \`a la seconde d\'emonstration du th\'eor\`eme 1.}

La fonction $G(x,y)=\int_{[0,1]}y(s)\,F'(x)(ds)$ est continue de $W\times W$ dans $\RR$ et v\'erifie $|G(x,y)|\leq \|y\|^2\sup_{x}\|F'(x)\|^2$. Le
 th\'eor\`eme 2 \'etendu aux variables bidimensionnelles, ce qui est sans difficult\'e,
s'applique et donne par le lemme 2
$$\EE\Gamma[F(X_n)]=\EE\hat{\EE}[((F(X_n))^\#)^2]\rightarrow\EE[\Gamma_0[F(B)]].\hspace{3cm}{\mbox{C.Q.F.D.}}$$

Par les propri\'et\'es de la convergence \'etroite vis \`a vis des fonctions Riemann-int\'egrables, cette seconde 
d\'emonstration donne \'egalement :\\

\noindent{\bf Corollaire 1.} {\it Soit $\Phi(x,y)$ une fonction de $W\times \hat{W}$ dans $\RR$ continue hors d'un
 n\'egligeable pour $m\times \mu$ o\`u $m$ est 
la loi de $B$ et $\mu$ celle de $B^\#$, telle que
$$|\Phi(x,y)|\leq K_1\|x\|^2+K_2\|y\|^2+K_3$$ alors $\EE\hat{\EE}\Phi(X_n,X_n^\#)\rightarrow\EE\hat{\EE}\Phi(B,B^\#)$.}\\

\noindent{\it Application.} Supposons ici pour simplifier les notations que $\sigma^2=c=1$ de sorte que $B^\#=\hat{B}$. La norme
 uniforme $N(w)=\|w\|$ qui est continue
et lipschitzienne appartient \`a $\DD_0$ (cf. [7] avec la m\'ethode Feyel-La Pradelle [9], ou [12] p.90), de m\^eme la fonctionnelle
 $M(w)=\sup_{t}w(t)$.

D'apr\`es les r\'esultats de Nualart et Vives [11] les op\'erateur $(.)^\#$ et $\Gamma_0$ sont calculables sur ces fonctionnelles :

i) $M^\#(w,\hat{w})=\hat{B}_\Sigma=\hat{w}(\Sigma(w))$ o\`u $\Sigma=\inf\{t:B(t)=\sup_sB(s)\}$
d'o\`u il r\'esulte que $\Gamma_0[M]=\Sigma$.

ii) $N^\#(w,\hat{w})={\mbox{sign}}(B_{\cal T})\hat{B}_{\cal T}$ o\`u ${\cal T}=\inf\{t:|B(t)|=\sup_s|B(s)|\}$
 d'o\`u il r\'esulte que $\Gamma_0[N]={\cal T}$.\\

L'ensemble des trajectoires browniennes qui atteignent plusieurs fois leur maxi\-mum est n\'egligeable et hors de cet 
ensemble, il n'est pas difficile de voir que l'application
$w\mapsto \Sigma(w)$ est continue, il r\'esulte alors du corollaire 1 que lorsque $n\uparrow\infty$ :
$$\EE\Gamma[\sup_tX_n(t)]=\EE\Gamma[\frac{1}{\sqrt{n}}\max_{1\leq k\leq n}S_k]\rightarrow\EE\hat{\EE}[M^{\#2}]=\EE[\Sigma]$$
et de m\^eme
$$\EE\Gamma[\|X_n(t)\|_\infty]=\EE\Gamma[\frac{1}{\sqrt{n}}\max_{1\leq k\leq n}|S_k|]\rightarrow\EE\hat{\EE}[N^{\#2}]=\EE[{\cal T}].$$

Ainsi nous pouvons \'enoncer avec les notations ci-dessus :\\

\noindent{\bf Proposition 2.} {\it Soit $F:\RR^2\mapsto \RR$ de classe $C^1\cap Lip$, alors d'une part
$$
\EE[F^2(\frac{1}{\sqrt{n}}\max_{1\leq k\leq n}S_k,\frac{1}{\sqrt{n}}\max_{1\leq k\leq n}|S_k|)]\rightarrow\EE[F^2(M,N)]$$
d'autre part
$$
\begin{array}{l}
\EE\Gamma[F(\frac{1}{\sqrt{n}}\max_{1\leq k\leq n}S_k,\frac{1}{\sqrt{n}}\max_{1\leq k\leq n}|S_k|)]\quad\quad\\
\quad\quad\rightarrow\EE[F_1^{\prime 2}(M,N){\cal T}]+2\EE[F'_1(M,N)F'_2(M,N)\,{\cal T}\wedge\Sigma]+\EE[F_2^{\prime 2}(M,N)\Sigma].
\end{array}
$$}

\noindent{\it Ind\'ependance asymptotique.} Lorsque $n$ augmente ind\'efiniment, la prise en compte d'un nombre croissant de $U_n$ fait que le processus $X_n$
se comporte comme un processus ind\'ependant d'une variable fix\'ee \`a l'avance. Le r\'esultat suivant est classique :

\noindent{\bf Proposition 3.}{\it  Soit $Y$ une v. a. d\'efinie sur $(\Omega, {\cal A}, \PP)$ int\'egrable. Si $\varphi\in C_b(W,\RR)$ on a
$$\EE[Y\varphi(X_n)]\rightarrow\EE[Y]\EE[\varphi(B)].$$}

Dans le m\^eme esprit, gr\^ace aux th\'eor\`emes 1 et 2  on a 

\noindent{\bf Th\'eor\`eme 3.} {\it  Si $Z\in L^\infty$, $Y\in\DD\cap L^\infty$, $\psi\in C^1\cap Lip$ et born\'ee, alors
$$\EE[Z\Gamma[Y\psi(X_n)]]\rightarrow\EE[Z\Gamma[Y]]\EE[(\psi(B))^2]+\EE[ZY^2]\EE[\Gamma_0[\psi(B)]].$$}

Observons que si $Z\geq 0$ et $\EE Z=1$ sous la probabilit\'e ${\bf Q}=Z.\PP$ les variables $U_n$ ne sont plus n\'ecessairement ind\'ependantes, et
de plus la forme ${\bf Q}[Z\Gamma[.]]$ n'est plus n\'ecessairement fermable. Le r\'esultat pr\'ec\'edent est donc une extension
 stricte du th\'eor\`eme 1.\\

{\it Remarque finale}. Terminons par quelques mots sur le r\'esultat principal lui-m\^{e}me. Supposons que les $U_n$ soient simul\'ees
 par une m\'ethode
de Monte Carlo avec une certaine pr\'ecision, de telle fa\c{c}on que l'hypoth\`ese d'ind\'ependance et de stationarit\'e des variables et
 de leurs erreurs puisse \^{e}tre consid\'er\'ee comme accep\-table. Contrairement ˆ certains th\'eor\`emes limites comme la loi des grands nombres 
qui effacent les erreurs (cf [4]), la normalisation faite pour la convergence en loi vers le brownien ne conduit sur celui-ci ni \`a une erreur nulle ni \`a
 une erreur infinie mais \`a l'erreur d'Ornstein-Uhlenbeck. Que ce soit cette structure d'erreur qu'on obtienne se con\c{c}oit bien car, 
d'apr\`es la formule de Mehler
(cf [12] p49 et [2] p116 \S 2.5.9), l'erreur qu'elle d\'ecrit est transversale et stationnaire. Nous voyons donc que pour obtenir d'autres structures d'erreur
sur l'espace de Wiener, telles que des structures de Mehler g\'en\'eralis\'ees (cf [2] p113 \S 2.5), il faut supposer que les erreurs sur les $U_n$ sont
correl\'ees.

\begin{list}{}
{\setlength{\itemsep}{0cm}\setlength{\leftmargin}{0.5cm}\setlength{\parsep}{0cm}\setlength{\listparindent}{-0.5cm}}
  \item\begin{center}
{\small REFERENCES}
\end{center}\vspace{0.4cm}

 [1] {\sc Billingsley, P.} {\it Convergence of Probability measures}, Wiley 1968.

[2] {\sc Bouleau, N.} {\it Error Calculus for Finance ansd Physics, the Language of Dirichlet Forms}, De Gruyter, 2003.

[3] {\sc Bouleau, N.} ``Algunes consideracions sobre llenguatges axiomatitzats amb eines d'extensi\'o: un enfocament en la teoria de la probabilitat i 
el c\`alcul d'errors amb formes de Dirichlet" {\it Bull. de la Soc. Catalana de Matematiques} Vol. 18, n\'um. 2, 2003, 25-37.

 [4] {\sc Bouleau, N.} ``Calcul d'erreur complet lipschitzien et formes de Dirichlet", {\it J. Math. pures et
 appl.} 80, 9, 961-976, 2001

[5] {\sc Bouleau, N.} ``Error calculus and path sensitivity in Financial models", {\it Mathematical Finance} vol 13/1, jan 2003, 115-134.

[6] {\sc Bouleau, N. } et {\sc Chorro, Chr.} ``Error structures and parameter estimation" {\it C. R. Acad. Sci. Paris }s\'er I 338, (2004) 305-310.

 [7] {\sc Bouleau, N.,} et {\sc Hirsch, F.}  {\it Dirichlet forms and analysis on Wiener space,} De Gruyter, 1991.

[8] {\sc Donsker, M.} ``An invariance principle for certain probability limit theorems" {\it Mem. Amer. Math. Soc.} no 6 ,1951.

 [9] {\sc Feyel, D.,} et {\sc la Pradelle, A. de} : ``Espaces de Sobolev Gaussiens", {\it Ann. Inst. Fourier}, 39-4, 875-908,  1989.

 [10] {\sc Neveu, J.} {\it Martingales \`a temps discret}, Masson, 1972.

[11] {\sc Nualart, D.} et {\sc Vives, J.} ``Continuit\'e de la loi du maximum d'un processus continu" {\it C. R. Acd. Sci. Paris} s\'er I, 307, 349-354, 1988.

[12]{\sc Nualart, N.} : {\it 
The Malliavin calculus and related topics}. Springer, 1995.
\end{list}

\end{document}